\documentclass[12pt,english]{article}
\usepackage[margin=0.75in]{geometry}
\usepackage[utf8]{inputenc}
\usepackage[english]{babel}
\usepackage{amsthm}
\usepackage{graphicx}
\usepackage{url}
\usepackage{hyperref}
\usepackage{amsmath}
\usepackage{listings}
\usepackage{amsfonts}
\usepackage{xcolor}

\begin{document}

\title{On the series expansion of a square-free zeta series}

\author{Artur Kawalec}

\date{}
\maketitle

\begin{abstract}
In this article, we develop a square-free zeta series associated with the M\"obius function into a power series, and prove a Stieltjes like formula for these expansion coefficients. We also investigate another analytical continuation of these series and develop a formula for $\zeta(\tfrac{1}{2})$ in terms of the M\"obius function, and in the last part, we explore an alternating series version of these results.
\end{abstract}

\section{Introduction}
In the recent work of Wolf [7] who investigated the divergence of a certain zeta series involving the M\"obius function $\mu(n)$ [3, p. 304], (thus essentially running over square-free integers) as such

\begin{equation}\label{eq:1}
\sum_{n=1}^{x}\frac{|\mu(n)|}{n}=O(\log x)
\end{equation}
(as $x\to \infty$) and determined a more refined estimate

\begin{equation}\label{eq:1}
\sum_{n=1}^{x}\frac{|\mu(n)|}{n}=\gamma^M+\frac{6}{\pi^2}\log x +O(\frac{1}{x}),
\end{equation}
where a new constant

\begin{equation}\label{eq:1}
\begin{aligned}
\gamma^M &=\lim_{x\to\infty} \left(\sum_{n=1}^{x}\frac{|\mu(n)|}{n}-\frac{6}{\pi^2}\log x\right)\\
&=1.04389451571193829740\ldots ,
\end{aligned}
\end{equation}
can be extracted in the limit in a much the same way as the classical Euler-Mascheroni constant

\begin{equation}\label{eq:1}
\begin{aligned}
\gamma &=\lim_{x\to\infty} \left(\sum_{n=1}^{x}\frac{1}{n}-\log x \right) \\
&=0.57721566490153286060\ldots
\end{aligned}
\end{equation}

Recalling that the Riemann zeta function is defined by the simplest Dirichlet series

\begin{equation}\label{eq:1}
\zeta(s)=\sum_{n=1}^{\infty}\frac{1}{n^s}
\end{equation}
is absolutely convergent for $\Re(s)>1$, and admits the Laurent series expansion about $s=1$ is given by

\begin{equation}\label{eq:1}
\zeta(s)=\frac{1}{s-1}+\sum_{n=0}^{\infty}\gamma_n \frac{(-1)^n(s-1)^n}{n!},
\end{equation}
which analytically extends $\zeta$ to $\mathbb{C}\backslash 1$ with residue $1$ and an only pole at $s=1$.  With the Euler's constant ($\gamma_0=\gamma$) appearing as the $0^{th}$ order coefficient in the series, and consequently, the next higher order Stieltjes constants $\gamma_n$ [4, p. 561] can be similarly generated by the well-known formula
\begin{equation}\label{eq:1}
\gamma_n=\lim_{x\to\infty} \Bigg\{\sum_{k=1}^{x}\frac{\log^n(k)}{k}-\frac{\log^{n+1}(x)}{n+1}\Bigg\}.
\end{equation}
And in a similar fashion, a certain fraction of the zetas

\begin{equation}\label{eq:1}
\frac{\zeta(s)}{\zeta(2s)}=\sum_{n=1}^{\infty}\frac{|\mu(n)|}{n^s}
\end{equation}
involves the M\"obius function as in (1) [6, p.5], and based on some previous analysis its Laurent series expansion is given by

\begin{equation}\label{eq:1}
\frac{\zeta(s)}{\zeta(2s)}=\frac{6}{\pi^2(s-1)}+\sum_{n=0}^{\infty}\gamma_n^M \frac{(-1)^n(s-1)^n}{n!},
\end{equation}
also with a pole at $s=1$ and residue $\frac{6}{\pi^2}$, and the proposed formula by Wolf for these expansion coefficients
\begin{equation}\label{eq:1}
\gamma^M_n =\lim_{x\to\infty}\Bigg\{\sum_{k=1}^{x}\frac{|\mu(k)|\log^n(k)}{k}-\frac{6}{\pi^2}\frac{\log^{n+1}(x)}{n+1}\Bigg\}
\end{equation}
is an analogue of the Stieltjes formula (7) [7]. We also note that the radius of convergence of (9) is only limited to $R=2$ due to the pole at the first trivial zero when $2s=-2$ for $s=-1$ from the center at $s=1$.

We now prove (10) by following exactly the proof in Bohman-Fr\"oberg [2][4, p. 561-562], but additionally inserting the M\"obius function. We begin the proof by introducing a self-canceling telescoping series

\begin{equation}\label{eq:1}
\sum_{k=1}^{\infty}(|\mu(k)|k^{1-s}-|\mu(k+1)|(k+1)^{1-s})=1
\end{equation}
then one has
\begin{equation}\label{eq:1}
(s-1)\frac{\zeta(s)}{\zeta(2s)}=1+\sum_{k=1}^{\infty}\left(|\mu(k+1)|(k+1)^{1-s}-|\mu(k)|k^{1-s}+(s-1)|\mu(k)|k^{-s}\right)
\end{equation}
which is still valid for $\Re(s)>1$, but since for $s=1$ the (lhs) limit is

\begin{equation}\label{eq:1}
\lim_{s\to 1}(s-1)\frac{\zeta(s)}{\zeta(2s)}=\frac{6}{\pi^2}
\end{equation}
is inconsistent with (rhs) of (12). So if we consider another re-scaled sequence instead

\begin{equation}\label{eq:1}
(s-1)\frac{\zeta(s)}{\zeta(2s)}=\frac{6}{\pi^2}+\sum_{k=1}^{\infty}\left(\frac{6}{\pi^2}|\mu(k+1)|(k+1)^{1-s}-\frac{6}{\pi^2}|\mu(k)|k^{1-s}+(s-1)|\mu(k)|k^{-s}\right)
\end{equation}
such that the (lhs) equals (rhs) at $s=1$, and then proceeding with the exp-log expansion exactly about $s=1$, then one obtains

\begin{equation}\label{eq:1}
\begin{aligned}
(s-1)\frac{\zeta(s)}{\zeta(2s)}=& \frac{6}{\pi^2}+\sum_{k=1}^{\infty}\Big[\frac{6}{\pi^2}\exp(-|\mu(k+1)|(1-s)\log(k+1))+\\
&-\frac{6}{\pi^2}\exp(|\mu(k)|(1-s)\log(k))+(s-1)\frac{1}{k}\exp(-(s-1)|\mu(k)|\log(k))\Big]
\\
\\
=& \frac{6}{\pi^2}+\sum_{k=1}^{\infty}\Bigg[\frac{6}{\pi^2}\sum_{n=0}^{\infty}\frac{(-1)^n(s-1)^n}{n!}\Big[|\mu(k+1)|^n\log^n(k+1)-|\mu(k)|^n\log^n(k)\Big]+\\
&+(s-1)\frac{1}{k}\sum_{n=0}^{\infty}\frac{(-1)^n(s-1)^n|\mu(k)|^n\log^n(k)}{n!}\Bigg]
\end{aligned}
\end{equation}
and by collecting the $(s-1)$ terms (with $(s-1)^0=1$) and since $|\mu(k)|^n=|\mu(k)|$, we obtain

\begin{equation}\label{eq:1}
\frac{\zeta(s)}{\zeta(2s)}=\frac{6}{\pi^2(s-1)}+\sum_{n=0}^{\infty}\gamma^M_n \frac{(-1)^n(s-1)^n}{n!}
\end{equation}
where

\begin{equation}\label{eq:1}
\gamma^M_n =\sum_{k=1}^{\infty}\left[\frac{|\mu(k)|\log^n(k)}{k}-\frac{6}{\pi^2}\frac{\log^{n+1}(k+1)-\log^{n+1}(k)}{k}\right]
\end{equation}
and last term as self-cancels again leaving only the $n+1$ term leading to the final limit formula

\begin{equation}\label{eq:1}
\gamma^M_n =\lim_{x\to\infty}\Bigg\{\sum_{k=1}^{x}\frac{|\mu(k)|\log^n(k)}{k}-\frac{6}{\pi^2}\frac{\log^{n+1}(x)}{n+1}\Bigg\}.
\end{equation}

We next consider an alternative representation for these coefficients. If we take the Laurent series expansion about $s=1$ of

\begin{equation}\label{eq:1}
\zeta(s)=\frac{1}{s-1}+\gamma+O(|s-1|)
\end{equation}
to the $0^{th}$ order, and another expansion of

\begin{equation}\label{eq:1}
\frac{1}{\zeta(2s)}=\frac{1}{\zeta(2)}-2\frac{\zeta^{\prime}(2)}{\zeta(2)^2}(s-1)+O(|s-1|^2)
\end{equation}
to the $1^{st}$ order, then by multiplying them together yields

\begin{equation}\label{eq:1}
\frac{\zeta(s)}{\zeta(2s)}=\frac{6}{\pi^2(s-1)}+\frac{\gamma}{\zeta(2)}-2\frac{\zeta^{\prime}(2)}{\zeta(2)^2}+O(|s-1|).
\end{equation}
As a result, the $0^{th}$ order coefficient in (21) is collected in the limit $s\to 1$ as

\begin{equation}\label{eq:1}
\gamma^M =\frac{6\gamma}{\pi^2}-\frac{72\zeta^{\prime}(2)}{\pi^4}
\end{equation}
gives a better closed-form formula in terms of other known constants. And the value for the zeta derivative $\zeta'(2)$ can be computed directly from (5) as

\begin{equation}\label{eq:1}
\zeta^{\prime}(2)=\sum_{k=1}^{\infty}\frac{\log(k)}{k^2}.
\end{equation}

In Table 1, we compute the first $10$ higher order $\gamma^M_n$ coefficients to high-precision ($30$ decimal places) by the $n^{th}$ differentiation of the (lhs) of (16) in the limit as

\begin{equation}\label{eq:1}
\gamma^M_n =(-1)^n\frac{d^n}{ds^n}\Bigg[\frac{\zeta(s)}{\zeta(2s)}-\frac{6}{\pi^2(s-1)}\Bigg]\Bigr\rvert_{s\to 1}
\end{equation}
which is easily possible to do in software packages such as Mathematica or Pari/GP [5][8]. The error term for $\gamma^M$ in (1) is $O(\frac{1}{x})$, and this would roughly imply that in order to get $30$ decimal places using (3), one needs to compute the limit with $x$ on the order of $\sim 10^{30}$, which is so large that even on modern computer workstations is still completely unfeasible to do in reasonable time (but perhaps on a supercomputer it could run faster), but one can still easily compute these constants to thousands of digits by differentiation in (24).

\begin{table}[hbt!]
\caption{A high-precision computation of expansion coefficients} 
\centering 
\begin{tabular}{| c | c |} 
\hline 
$n$ & $\gamma_n^M$ \\ [0.5ex] 
\hline 
$0$ & 1.043894515711938297404563438509   \\
\hline
$1$ &0.236152886477122974860578286060  \\
\hline
$2$ & 0.319384120408014249249465207074  \\
\hline
$3$ &0.501294458741649566645935631332   \\
\hline
$4$ & 1.010739722784850417039579626049   \\
\hline
$5$ &2.544030257932552280334481508980   \\
\hline
$6$ & 7.666100995112318690725728704276  \\
\hline
$7$ &26.88797470534219199661349019865  \\
\hline
$8$ & 107.6566910334506652692812639473  \\
\hline
$9$ &484.6934692784684121614213582581  \\
\hline
$10$& 2424.080089640181055133479838894  \\
\hline
\end{tabular}
\label{table:nonlin} 
\end{table}

\newpage

\section{On analytical continuation of the Dirichlet series}
We saw earlier that the square-free Dirichlet series

\begin{equation}\label{eq:1}
\frac{\zeta(s)}{\zeta(2s)}=\sum_{n=1}^{\infty}\frac{|\mu(n)|}{n^s}
\end{equation}
is convergent for $\Re(s)>1$. In the previous Section, the telescoping series (12) is also convergent for $\Re(s)>1$, but with the introduction of a scaling constant we propose that (14) is actually valid for $\Re(s)>\frac{1}{4}$ (except at $s=1$) and assuming (RH). We reformulate the expression and bring over $(s-1)$ to the (rhs) we have

\begin{equation}\label{eq:1}
\frac{\zeta(s)}{\zeta(2s)}=\frac{6}{\pi^2(s-1)}+\frac{1}{s-1}\sum_{k=1}^{\infty}\left[\frac{6}{\pi^2}\left(|\mu(k+1)|(k+1)^{1-s}-|\mu(k)|k^{1-s}\right)+(s-1)|\mu(k)|k^{-s}\right]
\end{equation}
and in another form

\begin{equation}\label{eq:1}
\frac{\zeta(s)}{\zeta(2s)}=\lim_{x\to\infty}\Bigg\{\sum_{k=1}^{x}\frac{|\mu(k)|}{k^s}+\frac{6}{\pi^2(s-1)}\left[1+\sum_{k=1}^{x}\left(|\mu(k+1)|(k+1)^{1-s}-|\mu(k)|k^{1-s}\right)\right]\Bigg\}
\end{equation}
and then self-canceling the telescoping series to last term we get

\begin{equation}\label{eq:1}
\frac{\zeta(s)}{\zeta(2s)}=\lim_{x\to\infty}\Bigg\{\sum_{n=1}^{x}\frac{|\mu(n)|}{n^s}-\frac{6}{\pi^2(1-s)}x^{1-s}\Bigg\}
\end{equation}
which we numerically find is valid $\Re(s)>\frac{1}{4}$ as stated. This limitation is due to the poles that arise at $\frac{1}{2}\rho$, where $\rho$ is the non-trivial root of $\zeta$. As evidenced by numerical computations, for example, for $s=0.8$ we compute the (lhs) of (28) to $-1.9413794172\ldots$, while the (rhs) with $x=10^8$ is $-1.941\underline{3}8634\ldots$, an agreement to $4$ digits. And if $s=0.5$ then the (lhs) of (28) is $0$, while the (rhs) with $x=10^8$ is computed as $-0.00173997\ldots$, where we see it is going to zero. And we checked many more random points and see that this formula is clearly converging to right values in the new domain $\Re(s)>\frac{1}{4}$, and as shown by the plot in Fig. 1, where we compare the (lhs) and (rhs) of equation (28) for $x=10^8$, and observe a deviation starts to happen near $\frac{1}{4}$. This means that all non-trivial roots are also roots of the (rhs) of (28).

\begin{figure}[h]
  \centering
  \includegraphics[width=150mm]{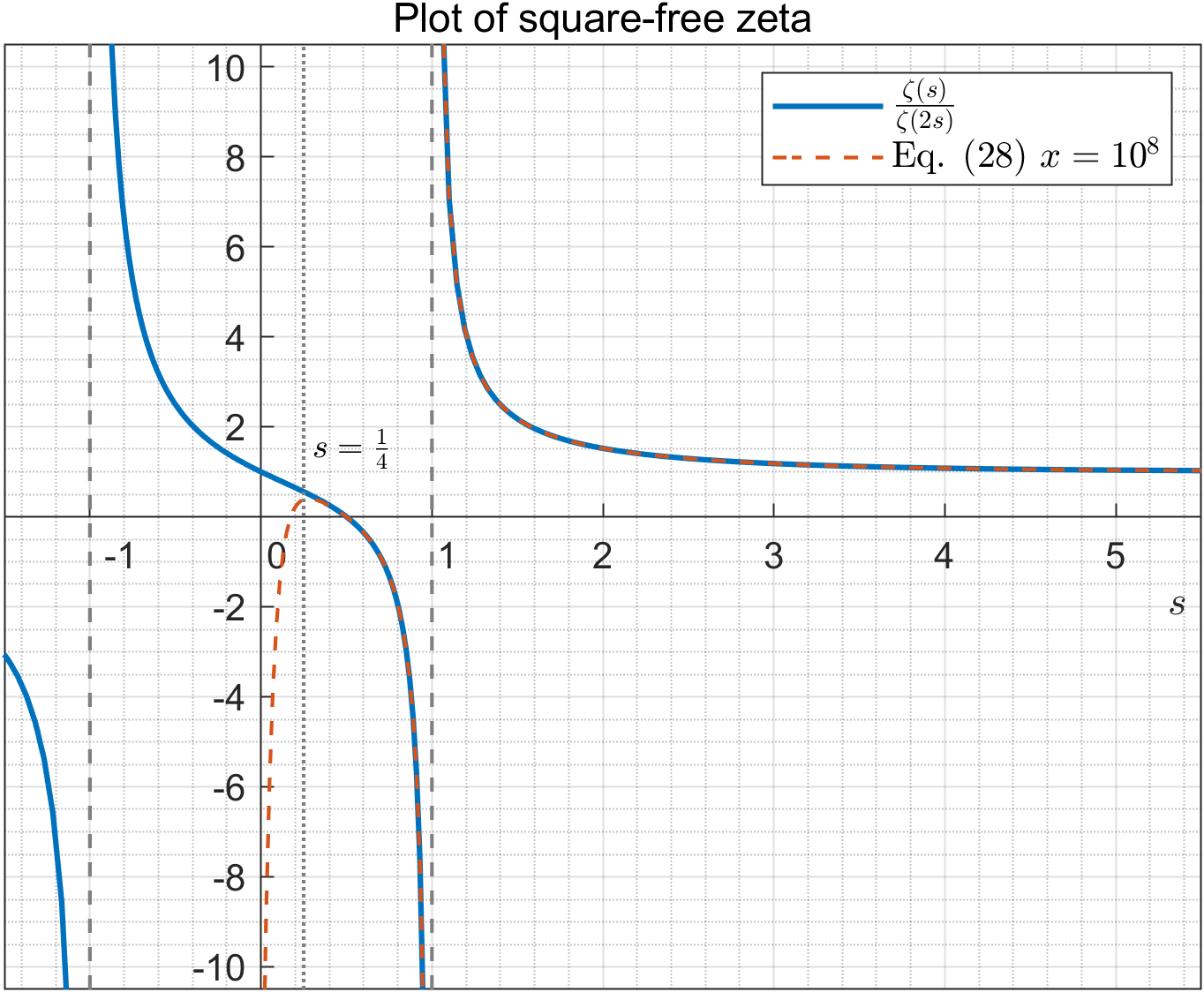}\\
  \caption{A plot of equation (28) for $x=10^8$ showing deviation near $s=\frac{1}{4}$}\label{1}
\end{figure}

Equation (28) can give many more interesting results. As previously seen,  $s=\frac{1}{2}$ is a zero of (28), and so we obtain a special case

\begin{equation}\label{eq:1}
\sum_{n=1}^{x}\frac{|\mu(n)|}{\sqrt{n}}\sim \frac{12}{\pi^2}\sqrt{x}
\end{equation}
as $x\to \infty$, or expressing in another way

\begin{equation}\label{eq:1}
\lim_{x\to\infty}\frac{1}{\sqrt{x}}\sum_{n=1}^{x}\frac{|\mu(n)|}{\sqrt{n}}= \frac{12}{\pi^2}.
\end{equation}

We also consider another expansion about $s=\frac{1}{2}$ of

\begin{equation}\label{eq:1}
\frac{\zeta(s)}{\zeta(2s)}=2\zeta(\frac{1}{2})(s-\frac{1}{2})+O(|s-\frac{1}{2}|^2)
\end{equation}
and seek to extract the first order zeta value by differentiation of the (rhs) of (28). We obtain the formula

\begin{equation}\label{eq:1}
\zeta(\frac{1}{2})=\lim_{x\to\infty}\Bigg\{-\frac{1}{2}\sum_{n=1}^{x}\frac{|\mu(n)|\log(n)}{\sqrt{n}}+\frac{6}{\pi^2}\sqrt{x}\left(-2+\log(x)\right)\Bigg\}
\end{equation}
which is converging to this zeta value. We next verify this formula numerically and find that it is a rather noisy and fluctuating series. In Fig. 2, we plot (32) as a function of $x$ from $x=10^0$ to $10^8$ in increments of $1$, on a logarithmic scale on the x-axis to capture all points.  We note how it is fluctuating about the value $\zeta(\frac{1}{2})=-1.4603545088\ldots$, with a gradual decrease in amplitude as $x$ increases. Another form can be simplified by plugging in (30) into (32) to have

\begin{equation}\label{eq:1}
\zeta(\frac{1}{2})=\lim_{x\to\infty}\Bigg\{-\frac{1}{2}\sum_{n=1}^{x}\frac{|\mu(n)|(\log(n)+2)}{\sqrt{n}}+\frac{6}{\pi^2}\sqrt{x}\log(x)\Bigg\}
\end{equation}
and when we test this formula we reproduce an almost the same plot as in Fig. 2.

\begin{figure}[h]
  \centering
  \includegraphics[width=160mm]{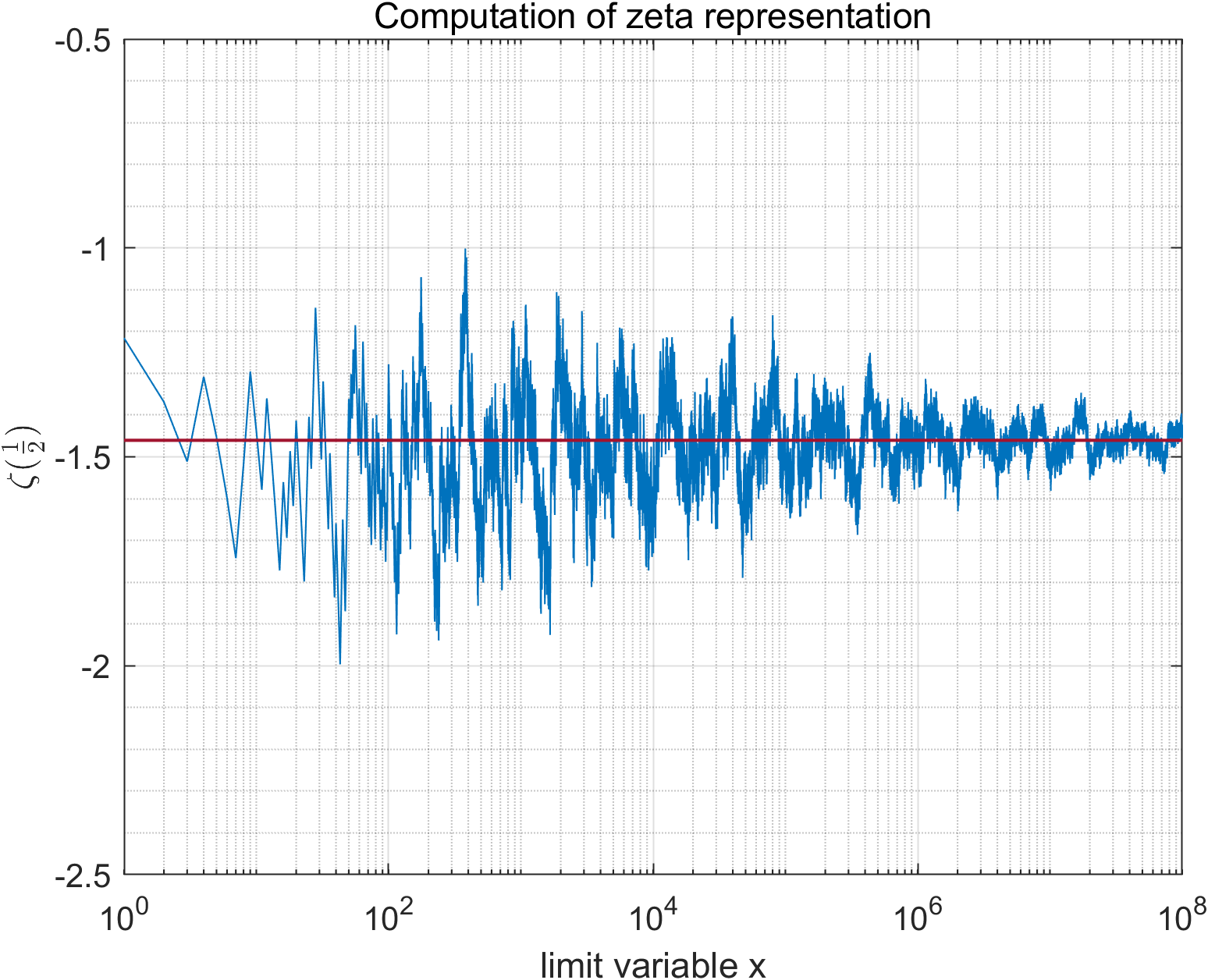}\\
  \caption{A plot of equation (32) for $\zeta(\frac{1}{2})$ as a function of limit variable $x$}\label{1}
\end{figure}

\newpage
\newpage

\section{On an alternating series representation}
There is a common alternating series representation for the zeta function

\begin{equation}\label{eq:1}
\zeta(s)=\frac{1}{1-2^{1-s}}\sum_{n=1}^{\infty}\frac{(-1)^{n+1}}{n^s}
\end{equation}
which is valid for $\Re(s)>0$ (but except at those points $t$ on the line $1+it$ where the leading factor becomes singular). A similar relation can be obtained for the square-free zeta series inspired by the MathOverflow post [9] by taking advantage of a certain connection to the Euler product

\begin{equation}\label{eq:1}
\frac{\zeta(s)}{\zeta(2s)}=\sum_{n=1}^{\infty}\frac{|\mu(n)|}{n^s}=\prod_{p}\left(1+\frac{1}{p^s}\right)
\end{equation}
valid for $\Re(s)>1$, where $p$ runs over all primes $p=\{2,3,5,7,\ldots\}$. Then, the application of the general Dirichlet series [3, p. 326] reads

\begin{equation}\label{eq:1}
\sum_{n=1}^{\infty}\frac{f(n)}{n^s}=\prod_{p}\sum_{k=0}^{\infty}\frac{f(p^k)}{p^{ks}},
\end{equation}
where $f(n)$ is a multiplicative function, meaning that for any positive integers $a$ and $b$ must satisfy the relation $f(1)=1$ and $f(a)f(b) = f(ab)$ whenever $a$ and $b$ are co-prime holds. It turns out that an alternating sign function $f(n)=(-1)^{n+1}$ is multiplicative, and hence $f(n)=(-1)^{n+1}|\mu(n)|$ is multiplicative, then we get the following transformation

\begin{equation}\label{eq:1}
\begin{aligned}
\sum_{n=1}^{\infty}(-1)^{n+1}\frac{|\mu(n)|}{n^s}&=\prod_{p}\sum_{k=0}^{\infty}\frac{(-1)^{p^k+1}|\mu(p^k)|}{p^{ks}}\\
&=\left(1-\frac{1}{2^s}\right)\prod_{p\geq 3}\left(1+\frac{1}{p^s}\right)\\
&=\left(1-\frac{1}{2^s}\right)\left(1+\frac{1}{2^s}\right)^{-1}\prod_{p\geq 2}\left(1+\frac{1}{p^s}\right) \\
&=\left(\frac{2^s-1}{2^s+1}\right)\frac{\zeta(s)}{\zeta(2s)}
\end{aligned}
\end{equation}
since $\mu(0)=1$, and $\mu(p)=-1$ and $\mu(p^k)=0$ for $k>1$. And this leads to the form analogue to the alternating zeta series (34) as

\begin{equation}\label{eq:1}
\frac{\zeta(s)}{\zeta(2s)}=\left(\frac{2^s+1}{2^s-1}\right)\sum_{n=1}^{\infty}(-1)^{n+1}\frac{|\mu(n)|}{n^s}.
\end{equation}
It is rare get such a coincidence, since this transformation depends on the existence of such a Euler product form. But unfortunately, the domain of convergence of (37) is only valid for $\Re(s)>1$ instead of $\Re(s)>0$ as in the alternating zeta (34). One reason is that at $s=1$ the summation in (34) is conditionally convergent as

\begin{equation}\label{eq:1}
\sum_{n=1}^{\infty}\frac{(-1)^{n+1}}{n}=\log(2)
\end{equation}
while the for square-free alternating zeta at $s=1$ diverges for the reason being by equating (8) with (34), where they only differ by a factor of $\frac{1}{3}$. And this leads to consider an estimate

\begin{equation}\label{eq:1}
\sum_{n=1}^{x}(-1)^{n+1}\frac{|\mu(n)|}{n}=\bar{\gamma}^M+\frac{2}{\pi^2}\log x +O(\frac{1}{x}),
\end{equation}
as $x\to \infty$, and an analogue constant for the alternating series

\begin{equation}\label{eq:1}
\begin{aligned}
\bar{\gamma}^M &=\lim_{x\to\infty} \left(\sum_{n=1}^{x}\frac{(-1)^{n+1}|\mu(n)|}{n}-\frac{2}{\pi^2}\log x\right)\\
&=0.53524615263113376955\ldots ,
\end{aligned}
\end{equation}
where the closed-form formula can be computed by

\begin{equation}\label{eq:1}
\bar{\gamma}^M =\frac{1}{3}\gamma^M  +\frac{8}{3\pi^2}\log(2)
\end{equation}
by the virtue of expansion of the factor
\begin{equation}\label{eq:1}
\frac{2^s-1}{2^s+1}=\frac{1}{3}+\frac{4}{9}\log(2)(s-1)+O(|s-1|^2)
\end{equation}
in conjunction with (19) to (22).

\section{Acknowledgement}
I would like to thank Professor Wolf for comments and suggestions for improving the paper.

\texttt{Email: art.kawalec@gmail.com}

\end{document}